\documentclass[11pt]{amsart}

\usepackage[english]{babel}

\usepackage{amsmath,amsthm, amssymb, bbm}
\usepackage{graphicx}
\usepackage{fullpage}

\newcommand{\R}{{{\mathbb {R}}}}

\newcommand{\Z}{{\mathcal Z}}
\newcommand{\z}{{\mathcal Z}}
\newcommand{\N}{{\mathcal N}}
\newcommand{\M}{{\mathcal M}}
\newcommand{\X}{{\mathcal X}}
\newcommand{\Y}{{\mathcal Y}}
\newcommand{\x}{{\mathcal X}}
\newcommand{\y}{{\mathcal Y}}
\newcommand{\n}{{\mathcal N}}

\newcommand{\U}{{\mathcal U}}

\newcommand{\G}{{\Gamma}}
\newcommand{\LL}{{\mathcal L}}

\newtheorem{Theorem}{Theorem}

\newtheorem {proposition}[Theorem]    {Proposition}

\let\oldmarginpar\marginpar
\renewcommand\marginpar[1]{\-\oldmarginpar[\raggedleft\scriptsize #1]%
{\raggedright\scriptsize #1}}

\begin{document}

\title{On the spatial Markov property of soups of unoriented and oriented loops}
\author{Wendelin Werner}
\address {
Department of Mathematics,
ETH Z\"urich, R\"amistr. 101,
8092 Z\"urich, Switzerland}
\email
{wendelin.werner@math.ethz.ch}

\begin {abstract}
We describe simple properties of some soups of {\em unoriented} Markov loops and of some soups of {\em oriented} Markov loops  that can be interpreted as 
a spatial Markov property of these loop-soups. This property of the latter soup is related to well-known features of the uniform spanning trees (such as Wilson's algorithm) 
while the Markov property of the former soup is related to the Gaussian Free Field and
 to identities used in the foundational papers of Symanzik, Nelson, and of Brydges, Fr\"ohlich and Spencer or Dynkin, or more recently by Le Jan. 
\end {abstract}

\maketitle

\section {Introduction}

Symanzik and then Nelson have pioneered the study of Euclidean field theory more than forty years ago \cite {Sy,Ne}. In their 
approach, measures on random paths and loops play an important role and led to further important developments such as in the work of Brydges, Fr\"ohlich and Spencer \cite {BFS} (see also Dynkin \cite {Dy,Dy2}). In all these papers, a gas of closed loops is used to represent partition functions and correlation structures of random fields.

The present note  will be in the same spirit, but the focus will be on this random gas of loops itself as the main object of interest, rather than viewing it as a combinatorial diagrammatic tool to evaluate quantities related to fields. We will in particular focus on the role of orientation of loops and describe a particular simple property of such random configurations of unoriented loops as well as for random configurations of  oriented loops. These properties are very directly related to the combinatorial features used in the aforementioned papers as well as to some features in the more recent study by Le Jan \cite {LJ}, who was also focusing more on properties of the occupation times of these soups. In particular, some of the observations 
in Sections 7 and 9 in \cite {LJ} can be viewed as describing some of the features that we will try to highlight. 

These gases of loops, or loop-soups (as they have been called in \cite {LW}) are a random Poissonian (i.e. non-interacting) collection of random unrooted loops in a domain, that can be associated naturally to a Markov process or a discrete-time Markov chain (see \cite {LJ} and the references therein).
 When one discovers the configurations of the loop-soup within a given sub-domain $U$ of the entire domain in which the soup is defined, one observes on the one hand loops that are entirely contained in $U$ (which form a loop-soup in $U$), and on the other hand, excursions in $U$ that are parts of loops that do not entirely stay in $U$. Note that different such excursions can belong to the same loop or not, depending on the configuration outside of $U$. The Markovian property that we shall discuss basically describes how to 
randomly complete the missing pieces into the loops i.e. it describes the conditional distribution of the loop-soup outside of $U$ when conditioning on these excursions  of the loop-soup in $U$. As we shall see, this takes a nice ``Markovian form'' in two special cases: 
\begin {itemize}
 \item 
 When one considers the loops to be oriented, and the intensity of the loop-soup to be the one that relates it to the partition function of uniform spanning trees ie. to the number of spanning trees (and to Wilson's algorithm \cite {Wi} to generate them uniformly at random, see e.g. \cite {Wi,La,LL,W1}). 
\item
 In the case where the chain is reversible, if one considers the loops to be unoriented, and chooses the intensity to be the one that relates the loop-soup to the Gaussian Free Field (for instance via their partition functions -- and in fact the occupation time of a continuous-time version of the loop-soup then corresponds exactly to the square of the GFF, see \cite {LJ}). 
\end {itemize}
In those two cases, the only relevant information in order to complete the excursions in $U$ into loops is the family of all endpoints of the excursions on $\partial U$, and not how these endpoints are connected by the excursions within $U$ (nor which excursion end-point is connected to which other by an excursion).  In other words, the trace of the discrete loop-soup inside $U$ and outside of $U$ are conditionally independent given their trace on $\partial U$ (more precisely, given their trace on the edges between $U$ and the complement of $U$).
 
\begin {figure}[ht!]
\begin {center}
\includegraphics[scale=0.6]{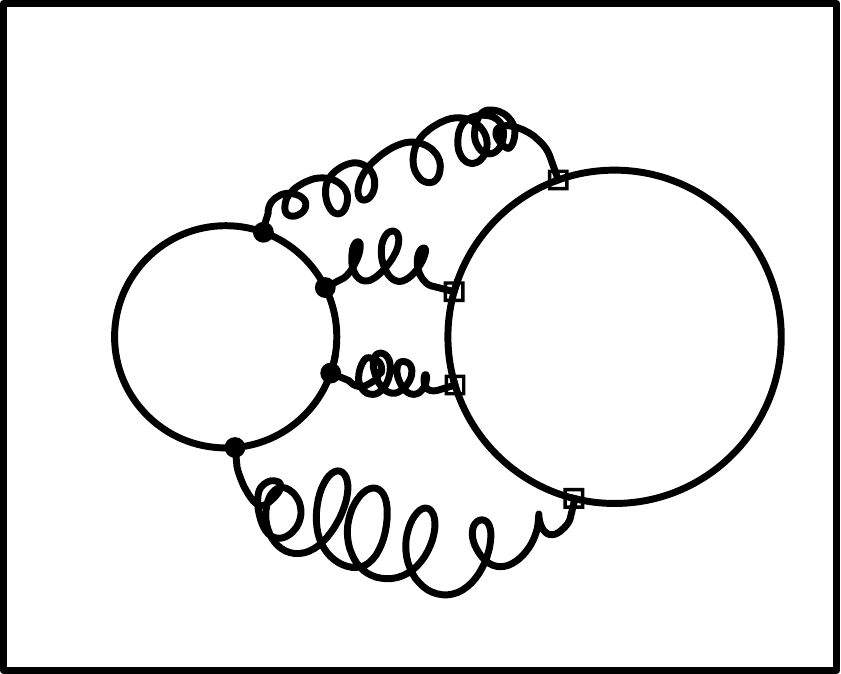}
\hskip 2mm
\includegraphics[scale=0.6]{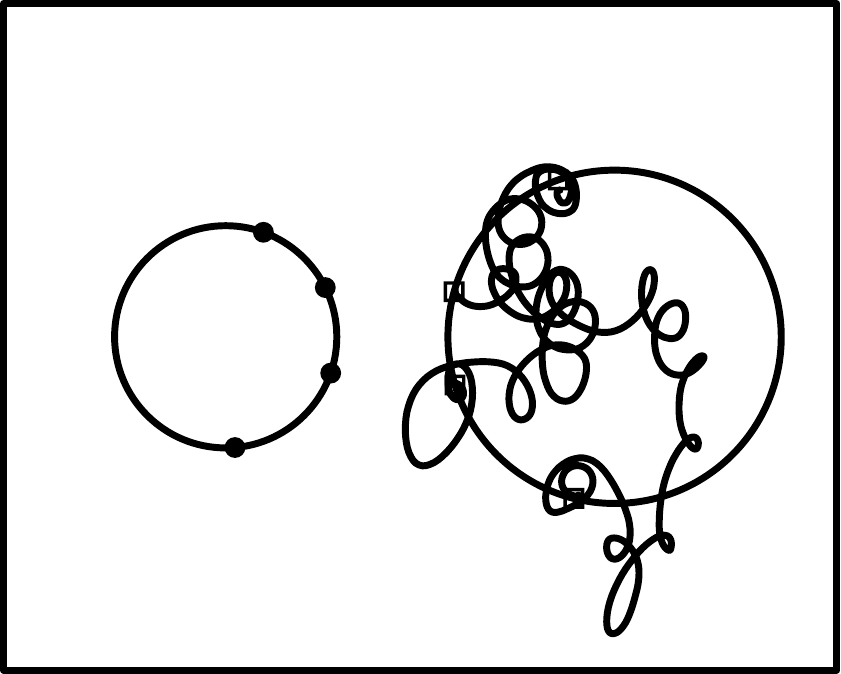}
\hskip 2mm
\includegraphics[scale=0.6]{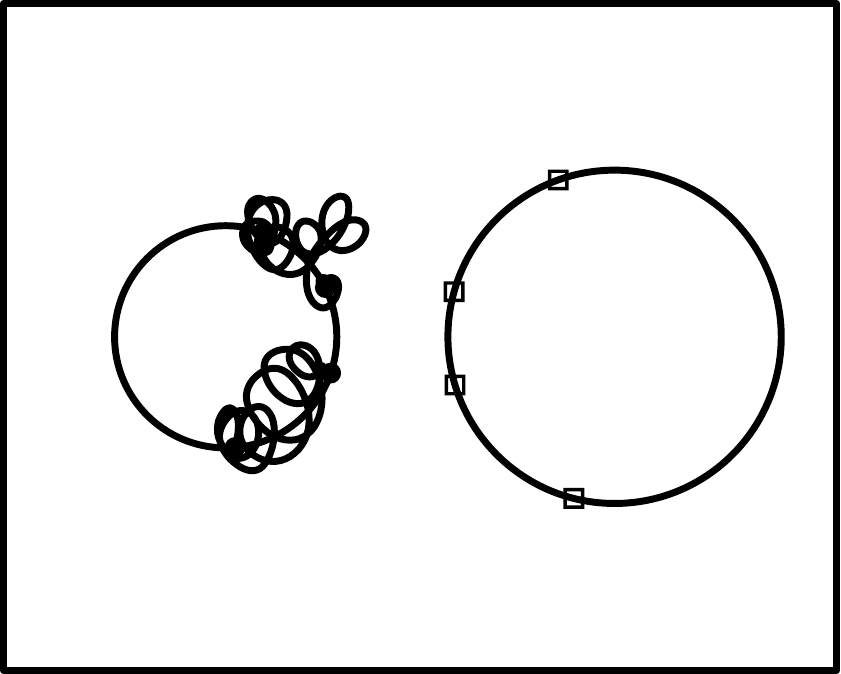}
\end{center}
\caption{Conditionally on the set of endpoints of crossings on each of the two circles, these three pictures, corresponding to different parts of the loops that touch both circles, are independent.}
\label {f2}
\end{figure}

Let us illustrate another instance of the spatial Markov property in an impressionistic and heuristic way via the following figures. We consider a loop-soup of unoriented loops in the inside of the rectangle, of well-chosen intensity (related to the partition function of the GFF). In this loop-soup, only finitely many loops do touch the two circles, and in each such loop, there are an even number of ``crossings'' from one circle to the other. The statement in the caption of Figure \ref {f2} is the type of result that we will derive. 
\begin{figure}[ht]
\begin{center}
\includegraphics[scale=0.6]{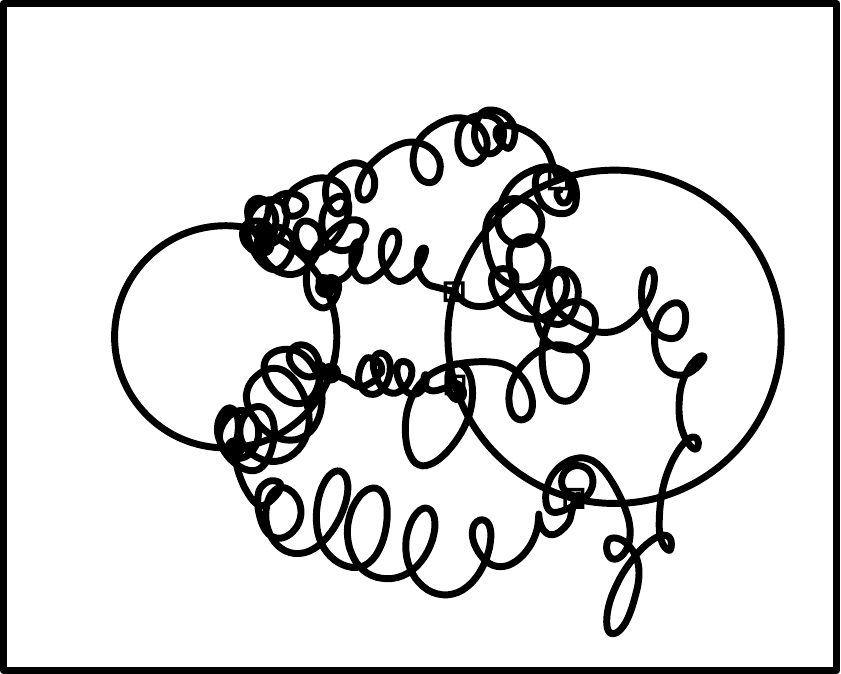}
\hskip 5mm
\includegraphics[scale=0.6]{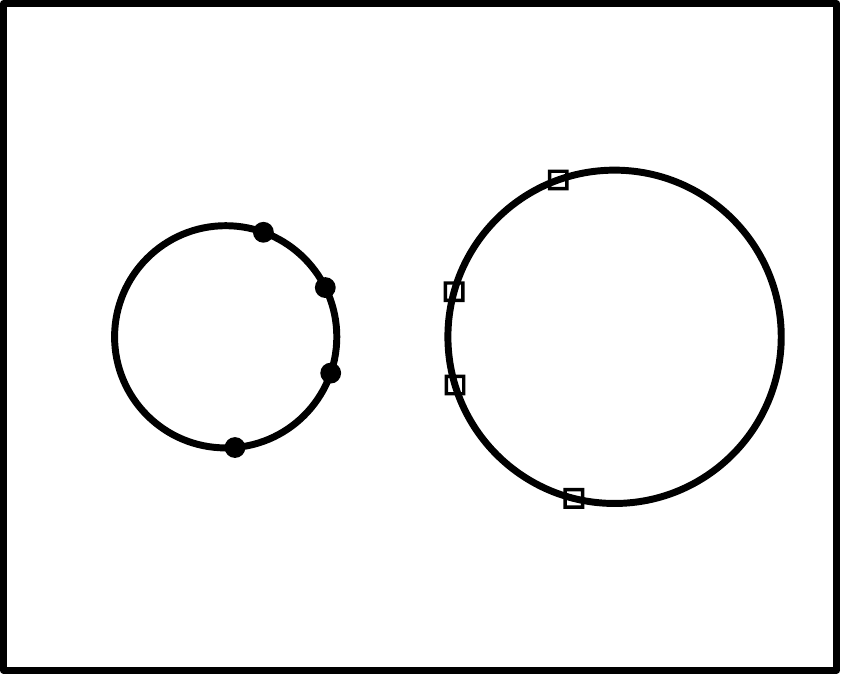}
\end {center}
\caption {The unoriented loop(s) in the soup that touch both circles, and the  endpoints of their (four in this case) crossings between the two circles}
\end {figure}

To conclude this introduction, let us briefly mention that of the motivations for the present work is to explore the relation between the natural ``Markovian'' structures emerging from the loop-soups with the theory of local sets for the discrete and continuous GFF, as defined by Schramm and Sheffield in \cite {SS}.

\section {Background and definitions}

In this section, we recall standard facts about Markov loops and loup-soups, make some elementary comments about the orientation/non-orientation of loops, and we define the natural measures on Markov bridges that we will need. 

\subsection{The measure on unrooted oriented loops}

Let us consider a discrete oriented graph $\G$, where each vertex $x$ has a finite number $d(x)$ of outgoing edges, so that it is possible to define simple random walk on $\G$ ($d(x)$ is however not necessarily the same for all $x$). Note that there could be ``several'' parallel edges from a vertex $x$ to a vertex $y$. Also, as opposed to the unoriented case, the could be an edge from $x$ to $y$ but no edge from $y$ to $x$.

We say that $l=(l_0, e_1, l_1, l_2, \ldots, l_{n-1}, e_{n-1}, l_n)$ for $n \ge 1$ is a rooted loop with $n = |l|$ steps in $\G$ if $l_0, l_1, l_{n-1}$ are sites of the graph, if $l_0 = l_n$ and if for all $i \in \{1, \ldots , n\}$, $e_i$ denotes an edge from $l_{i-1}$ to $l_{i}$ in the graph. Let us notice that in the case of parallel edges in the graph, the information about which oriented edges were used are part of the information contained in the loop.  

We can note that the probability $p(l)$ that a random walk starting from $l_0$ follows exactly this loop during its first $n$ steps is exactly $1/ \prod_{i=0}^{n-1} d(l_i)$. 
We define the measure $\rho$ on rooted loops by $\rho (l)= p(l) / n$. Note that this is not a probability measure (a loop $l$ might for instance contain another loop as its first steps if it visits $l_0$ several times before time $n$; furthermore, we sum over all possible starting points $l_0$ in the graph).

The quantity $\rho (l)$ remains unchanged if one changes the root of the loop (if one considers the loop $(l_i, e_{i+1},  l_{i+1}, \ldots, l_n, e_1, l_1, \ldots, l_i)$ instead of $l$), which  leads naturally to the definition of 
{\em an unrooted loop ${L}$} as an equivalence class of rooted loops, where two loops are equivalent as soon as they are obtained from one another by rerooting.
The measure $\mu$ on unrooted oriented loops is then the image of the measure $\rho$ under the mapping that maps each rooted oriented loop to its equivalence class of unrooted loops.
This is the loop-measure that has been used and studied extensively in recent years, in connection with loop-erased random walks, Gaussian Free Fields, Dynkin's isomorphism theorems and in the continuous two-dimensional (Brownian) setting, with conformal 
loop ensembles and SLE curves, see e.g. \cite {LJ,W1} and the references therein). 

In many cases, the number of different rooted loops in the same equivalence class of unrooted loops is the length $n(l)=|l|$ of the loop (one possible root per step on the loop). However, when a loop $l$ consist exactly of the concatenation 
of $J \ge 1$ copies of exactly the same loop, ie, $n=Jn_1$ and $l$ is exactly the concatenation of $J$ copies of $(l_0, \ldots , l_{n_1})$ (and $J=J(l)$ is the maximal such number -- note that this number is also invariant under rerooting of $l$ so that we can view it as a function of ${L}$), then the number of rooted loops that give rise to the same unrooted loop as $l$ is $n / J(L)$. Hence, the general formula for $\mu$ is 
$\mu ( {L} ) = p(l) / J(l)$, when $l$ is any loop in the equivalence class ${L}$. 

In the sequel, we will refer to loops $l=(l_0, e_1, \ldots, l_n)$ (or their equivalence class) such that $J(l)=1$ as single loops, 
and we say that the loop $l^k$ defined as the concatenation of $k$ copies of $l$ ie. as $(l_0, e_1, \ldots, l_{n-1}, l_0, e_1, \ldots, 
l_{n-1},  \ldots, l_{n-1}, e_n, l_0)$ with $J (l^k) = k$ is its $k$-fold multiple.

\subsection {The measure on unrooted unoriented loops}

In the previous subsection, the graph was oriented, and all our loops (rooted and unrooted) were oriented. 
Let us now consider an unoriented graph, where each vertex $x$ has a finite number $d(x)$ of outgoing edges (here a single edge from $x$ to $x$ would be counted twice, and we also allow parallel edges between two sites $x$ and $y$). 
Then the previous quantity $p(l)$ remains unchanged when one changes the orientation of the loop; indeed, if one defines the time-reversal   
$\hat l := (l_n, e_n, l_{n-1}, \ldots, l_1, e_1, l_0)$, then  $ p(l) = 1/ \prod_{i=1}^n d(l_i) = p (\hat l)$. 

We now define {an unrooted unoriented loop} as the equivalence class of oriented rooted loops, where two such loops are said to be equivalent as soon as they are obtained from one another by rerooting and possibly by time-reversal. Or alternatively, we say that an unrooted unoriented loop is the equivalence class of unrooted oriented loops, modulo time-reversal. 

We then define the  measure $\nu$ on unrooted unoriented loops to be the image of $\rho /2$ under the mapping that maps each rooted oriented loop onto to its equivalence class of unrooted and unoriented loops. The measure $\nu$ is of course just the unoriented projection of $\mu /2$.

When the time reversal $\hat l$ of a rooted oriented loop $l$ is not in the same unrooted oriented class of loops as $l$, then there will be twice more rooted oriented loops in the same class $\tilde   L$ of unoriented unrooted loops of $l$ than in its class $L$ of oriented unrooted loops, so that $\nu ( \tilde {L} ) =  \mu (L)$. 
It however can happen that $l$ and $\hat l$ define the same oriented unrooted loop ${L}$ (for instance when the loop $l$ is the concatenation of a loop with its time-reversal). In that case,  $\nu ( \tilde { L} ) =   \mu (L) / 2$. We define $\tilde J ( \tilde L) = J (L)$ or $2 J(L)$ depending on whether $L \not= \hat L$ or not, so that 
$\nu ( \tilde L) = p (l) / \tilde J (\tilde L)$ for all $\tilde L$.

\medbreak

All the previous definitions have also straightforward counterparts and generalizations for general  Markov processes (not necessarily random walks) -- the processes 
would need to be reversible for the unoriented loops --, and in continuous time and/or in continuous space. Note that as soon as one deals with continuous time, the multiplicity issues (raised by the fact that $J$ is not constant) do not exist. One fundamental example is of course the Brownian loop measure that gives rise to the loop-soup, as introduced in \cite {LW}. Other examples include the Brownian loops on cables systems associated to discrete graphs, as studied in \cite {Lupu1}. 

Since our  purpose here is to give an elementary presentation of the resampling property of loop-soups,  we have opted in the present paper to state and explain things in the most transparent settings (random walk loops
on regular graphs, where all points in $\G$ have the same number $g$ of outgoing edges -- which we will from now on assume --, and Brownian loops). The generalization of the 
proofs to continuous-time and discrete space Markov processes do not require any new idea.

\subsection {Loop-soups}

For a given graph, one can define  simple natural random objects out of the measures on loops. For each $\alpha >0$, one can define a Poisson point process of loops, with intensity given by $\alpha$ times the measure $\mu$ on loops. This is the loop-soup, as introduced in the Brownian setting in \cite {LW} and studied more recently in the discrete setting in \cite {LJ}. It is also the gas of loops that was already used in \cite {Sy,BFS}.

Of course, when one samples a soup of (unrooted) oriented loops according to the 
loop measure $\alpha \mu$, and one forgets about the orientation of the loops, one gets a soup of unrooted unoriented loops with intensity $2\alpha \nu$, and conversely, one can recover the former by choosing at random the orientation of each loop. 
In order to avoid confusions, we will use the letters $\alpha$ to denote the intensity of soups of oriented loops (i.e. with intensity measure $\alpha \mu$) and $c$ to denote the intensity of soups of unoriented loops (i.e. with intensity measure $c \nu$). The natural relation between $c$ and $\alpha$ is then $c= 2 \alpha$.

We will not recall all the properties of these loop-soups, but we would like to 
stress the following points:
\begin {itemize}
 \item 
 The  soup of oriented loops with intensity $\alpha=1$ is very closely related to uniform spanning trees. In particular, the loops in such a loop-soups correspond exactly to the family of loops that have been erased when performing Wilson's algorithm to sample a uniform spanning tree in $\G$. And in this context, it is somewhat more natural to consider oriented loops. 
\item
 The soup of unoriented loops with intensity $c=1$ is very closely related to the Gaussian Free Field in $\G$ and its square. In this context, because one looks only at the cumulated occupation times of the loops, it is in fact somewhat more natural to consider unoriented loops (as the orientation is not needed to define the occupation time measure).  
\end {itemize}
With this notation, the UST is related to $c=2$ and the GFF to $c=1$, and more generally, in two dimensions, in the conformal field theory language, the value of $c$ corresponds to the absolute value of the central charge of the corresponding models.
 
Suppose now that $L_1, \ldots, L_k$ are $k$ different oriented unrooted loops. Let $\U_1, \ldots, \U_k$ denote the respective number of occurrences of these loops in an unrooted loop-soup 
with intensity $\alpha \mu$. These are $k$ independent Poisson random variables with respective means $\alpha\mu (L_1), \ldots , 
\alpha\mu (L_k)$, so that 
$$ P ( \U_1 = u_1, \ldots ,  \U_{k}=u_k ) 
 = \prod_{j=1}^k  ((\alpha \mu ( L_j))^{u_j} e^{- \alpha \mu (L_j)} / u_j!).$$ 
In the special case where $\alpha=1$, the $\alpha^{u_j}$ terms disappear, and we get 
$$ \frac { P ( \U_1 = u_1, \ldots ,  \U_{k}=u_k )}{ P (\U_1 = \ldots = \U_k =0 ) }= \prod_{j=1}^k \frac {(p(L_j) / J(L_j))^{u_j} }{ u_j!}.$$

Similarly, if we are considering instead a loop-soup of unoriented loops with intensity $\nu$ (ie. for $c=1$), the very same formula holds, i.e. if $\tilde L_1, \ldots, \tilde L_k$ are $k$ different unoriented loops, and if $\tilde \U_1, \ldots, \tilde \U_k$ denote the respective number of occurrences of these loops in a soup of unrooted loops with intensity $\nu$, then 
$$ \frac { P ( \tilde \U_1 = u_1, \ldots ,  \tilde \U_{k}=u_k )}{ P (\tilde \U_1 = \ldots = \tilde \U_k =0 ) }= 
\prod_{j=1}^k \frac {(p(L_j) / \tilde J( \tilde L_j))^{u_j} }{ u_j!}.$$

\subsection {Random bridges}

Recall that in order to slightly simplify notations and some of our considerations, we are from now going to assume that (both 
in the oriented and in the unoriented cases), the graph $\G$ will be such that each site has the same number $g$ of outgoing edges. Note that this is not really a restriction, because it is for instance always possible starting from an unoriented graph $\G$ where each site $x$ has $d(x)$ outgoing edges, with $\sup_x d(x) \le g$, to add $(g - d(x))$ stationary edges from $x$ to $x$ to the graph, without changing the behavior of the random walks (and this leads to the natural way to extend the results to the case of graphs with non-constant degree). 
 
\medbreak

Let us first suppose that $\G$ is an oriented graph. 
Consider now a subgraph $D \subset \G$ and two points $x$ and $y$ in $D$.
We say that a bridge $b$ from $x$ to $y$ in $D$ is a finite nearest-neighbour path (keeping track of the oriented edges used) in $D$ starting at $x$ and finishing at $y$. We call $n(b)$ the length (number of jumps) of $b$. A bridge from $x$ to $x$ is allowed to have a zero length. 

Suppose now that the Green's function $G_D (x,y)$ is positive and finite. 
Recall that this is the mean number of visits at $y$ before exiting $D$, by a random walk starting at $x$. In other words, it is the sum over all bridges from $x$ to $y$ in $D$ of $g^{-n(b)}$. We can therefore define a probability measure  on bridges from $x$ to $y$ in $D$, that assigns a probability 
$ g^{-n(b)}  / G_D (x,y)$ to each bridge $b$. 

Suppose now that we are given $N$ points $x_1, \ldots , x_N$ and $N$ points $y_1, \dots , y_N$ in $D$. We  say that a the family of paths $b^1, \ldots , b^N$ is an ordered bridge in $D$ from $X=(x_1, \ldots, x_N)$ onto $Y=(y_1, \ldots , y_N)$ if each $b^j$ is a bridge from $x_j$ to $y_j$ in $D$. We also define $G_D ( X, Y) = \prod_{j=1}^N G_D (x_j, y_j)$ and when this quantity is not equal to zero nor infinite, we define the 
probability measure on ordered bridges from $X$ to $Y$ in $D$ to be obtained by taking $N$ independent bridges from $x_j$ to $y_j$ respectively.

An unordered bridge from $X$ to $Y$ is defined to be the knowledge of a permutation $s$ from $\{ 1 , \ldots N \}$ and of an ordered bridge from $X$ to 
$Y^s = (y_{s(1)}, \ldots, y_{s (N)})$. 
 We now define define a probability measure $B_{X,Y}^D$ on unordered bridges from $X$ to $Y$ in $D$ as follows: 
\begin {enumerate}
\item First sample a permutation $\sigma$ so that the probability of $\sigma=s$ is proportional to $G_D (X, Y^s)$. 
\item Then, conditionally on $\sigma$,  sample the ordered bridge from $X$ to $Y^\sigma$ according to the probability measure on ordered bridges in $D$ described above.
\end {enumerate}
For this to make sense, we need that for at least one $s$, $G_D ( X, Y^s) >0$. 
This procedure basically samples an unordered bridge from $X$ to $Y$ in such a way that the probability of a given unordered bridge is proportional to $g^{-K}$, where $K$ denote the sum of the length of the $N$ bridges that form the generalized bridge. 
Mind that in the present setting, when $y_2 = y_3$ say, we do count the same collection of $N$ bridges (corresponding to interchanging $y_2$ and $y_3$) twice in our partition function, because they correspond to different permutations.

\medbreak

Let us now suppose that the graph $\G$ is not oriented. In the previous definition,  each bridge has an implicit orientation (from $x$ to $y$). On the other hand, the image under time-reversal (i.e. consider $\hat b_j = b_{n-j}$) of the bridge probability from $x$ to $y$ in $D$  is exactly the bridge probability from $y$ to $x$ in $D$ (note that we use here the fact that $x$ and $y$ have the same number of outgoing edges $g$). One can therefore define the probability measure on unoriented bridges in $D$ joining $x$ and $y$ to be the law obtained by considering $B_{x,y}^D $ and then forgetting about the time-orientation.  

Suppose now that $Z=(z_1, \ldots, z_{2N})$ are $2N$ points in $D$. An unoriented $Z$-bridge is the knowledge of a pairing $t$ of $\{1, \ldots , 2N \}$    
  (this is a permutation that contains only cycles of length exactly $2$ -- and we say that $i$ and $t(i)$ are paired -- we will denote the $N$ pairs of $t$ by $(t^1_1, t^2_1), \ldots , (t_N^1, t_N^2)$ using some lexicographic rule), and of $N$ unoriented bridges joining the $N$ pairs $(z_{t_k^1}, z_{t_k^2})$ for $k \le N$.
  
   For each $Z$, we then define the measure $B_Z^D$ on unoriented unordered $Z$-bridges as follows: 
 \begin {enumerate}
  \item 
 We first sample a pairing $\tau$ in such a way that the probability of a given pairing $t$ is proportional to   $\prod_{k=1}^N G_D ( z_{t_k^1}, z_{t_k^2})$. 
  \item
  When $\tau=t$, we then sample an $N$ independent (unoriented) bridges in $D$ joining the two points of each of the $N$ pairs $(z_{t_k^1}, z_{t_k^2})$.
 \end {enumerate}
 Again, this only makes sense if for at least one pairing $t$, $\prod_k G_D (z_{t_k^1}, z_{t_k^2}) $ is positive. Then, the definition just means that we sample a $Z$-bridge in such a way that the probability of a given $Z$-bridge is just proportional to $g^{-K}$ where $K$ denote the sum of the length of the $N$ bridges that form this $Z$-bridge. 
 
 \medbreak
 
 These definitions of bridges can be trivially extended to the Brownian settings (both in $d$-dimensional space as well as on cable systems), provided that no two $z_j$'s coincide (in the unoriented bridges) and that no $x_i$ is equal to an $y_j$ (for the oriented bridges) so that the Green's functions involved are all finite. The only difference is that the distribution of an individual bridge from $x$ to $y$ is done in two steps: 
\begin {enumerate}
 \item 
First, sample the time-length $T$ of the Brownian bridge according to the probability measure $p_{D,t} (x,y) dt / G_D (x,y)$, where $p_{D,t} (x, y)$ is the density at
$y$ of the law of a Brownian motion at time $t$, starting from $x$ and killed upon exiting $D$.
\item Then, conditionally on $T$, sample a usual Brownian bridge from $x$ to $y$ and time-length $T$, conditioned to stay in $D$. 
\end {enumerate}

\section {Partial resampling of soups, and spatial Markov properties}

We now describe various instances of the partial resampling properties of loop-soups, and discuss some consequences. 

\subsection {Partial resampling of  soups of oriented loops at $\alpha=1$}
\label {ss2.1}

Let us suppose that $\G$ is an oriented graph of degree $g$ as before, and that $D \subset \G$ is a subgraph of $\G$ where the Green's function is finite. 
We are going to describe the resampling property of the soup of oriented  loops with intensity $\alpha=1$. 
Suppose that  $F_1$ and $F_2$ are two disjoint finite set of vertices in our graph. 
When one considers a loop-soup in $D$, then the number of loops in the loop-soup that do intersect both $F_1$ and $F_2$ is a Poisson random variable ${\mathcal M}={\mathcal M} (F_1, F_2)$ with finite mean equal 
to the $\mu$-mass of the set of loops that intersect both $F_1$ and $F_2$.
We denote the family of $\M$ loops that intersect both $F_1$ and $F_2$ by $\overline {\mathcal L}$ (the information in $\overline {\mathcal L}$ includes how many occurrences of any given oriented unrooted loop that intersects $F_1$ and $F_2$ there are). We will write $\overline {\mathcal L} = ( \LL_1, \ldots, \LL_\M)$, where the chosen order of the loops in the family follows some lexicographic (deterministic) rule, so that the information provided by $\overline \LL$ and $(\LL_1, \ldots, \LL_\M)$ are identical.

When $L$ is an unrooted loop that intersects $F_1$ and $F_2$, 
we can consider the finitely many portions of $L$ that are of the type $(a_0, e_1, a_1,  a_1, \ldots, a_k)$ with $a_0,  a_k \in F_2$, $\{ a_1, \ldots, a_{k-1} \} \cap F_2 = \emptyset$ and at least one of the $a_i$ is in $F_1$. In other words, these are the excursions of $L$ away from $F_2$ that do reach $F_1$. We allow $a_0=a_k$, or the excursion to be the entire loop (which happens if $L$ visits $F_2$ only once) and it can also happen that the same excursion occurs several times in the same loop. 

When we sample $\overline {\LL}$,  we call $\eta$ the collection of all excursions of its loops. We can again decide to order them in some lexicographic predetermined deterministic way, so that we can write  $\eta= ( \eta_1, \ldots , \eta_\N)$ (again, it is important that if a given piece appears several times in the loop-soup, then it appears several times in 
this list as well). 
Note that $\N \ge \M$ because each loop that intersects $F_1$ and $F_2$ contains at least one such excursion. The pieces $\eta_1, \ldots, \eta_\N$ might be part of $\N$ different loops (in which case $\N= \M$), but they could also be all parts of the same loop (in which case $\M=1$). Of course, the probability that $\N=\M= 0$ is also positive.  

Observe that one intuitive way to discover all these excursions is in fact to explore all the loops ``starting'' from their intersection points with $F_1$, in both the positive time-direction and the negative time-direction, until reaching $F_2$ in both directions. 

Each of the pieces $\eta_j$ are naturally oriented as parts of oriented loops, and we can define their respective starting points $\y_j$ and endpoints $\x_j$ (note that all these points are on $F_2$). The missing parts of the loops that the $\eta$'s are part of will therefore be bridges in the complement of $F_1$, that join each of the $\x_j$'s to a $\y_{\sigma (j)}$ for a permutation $\sigma$ ie. the missing part will be an unordered bridge $\beta$ from the vector $\X = (\X_1, \ldots, \X_\N)$ to the vector $\Y
= (\Y_1, \ldots, \Y_\N)$ in $D \setminus F_1$. 
Now, the resampling result in this case goes as follows: 

\begin {proposition}
\label {p1}
 The conditional distribution of $\beta$ given $\eta$ is exactly the unordered bridge measure $B_{\X,\Y}^{D \setminus F_1}$.
\end {proposition}

\begin{figure}[ht]
\begin{center}
\includegraphics[scale=0.5]{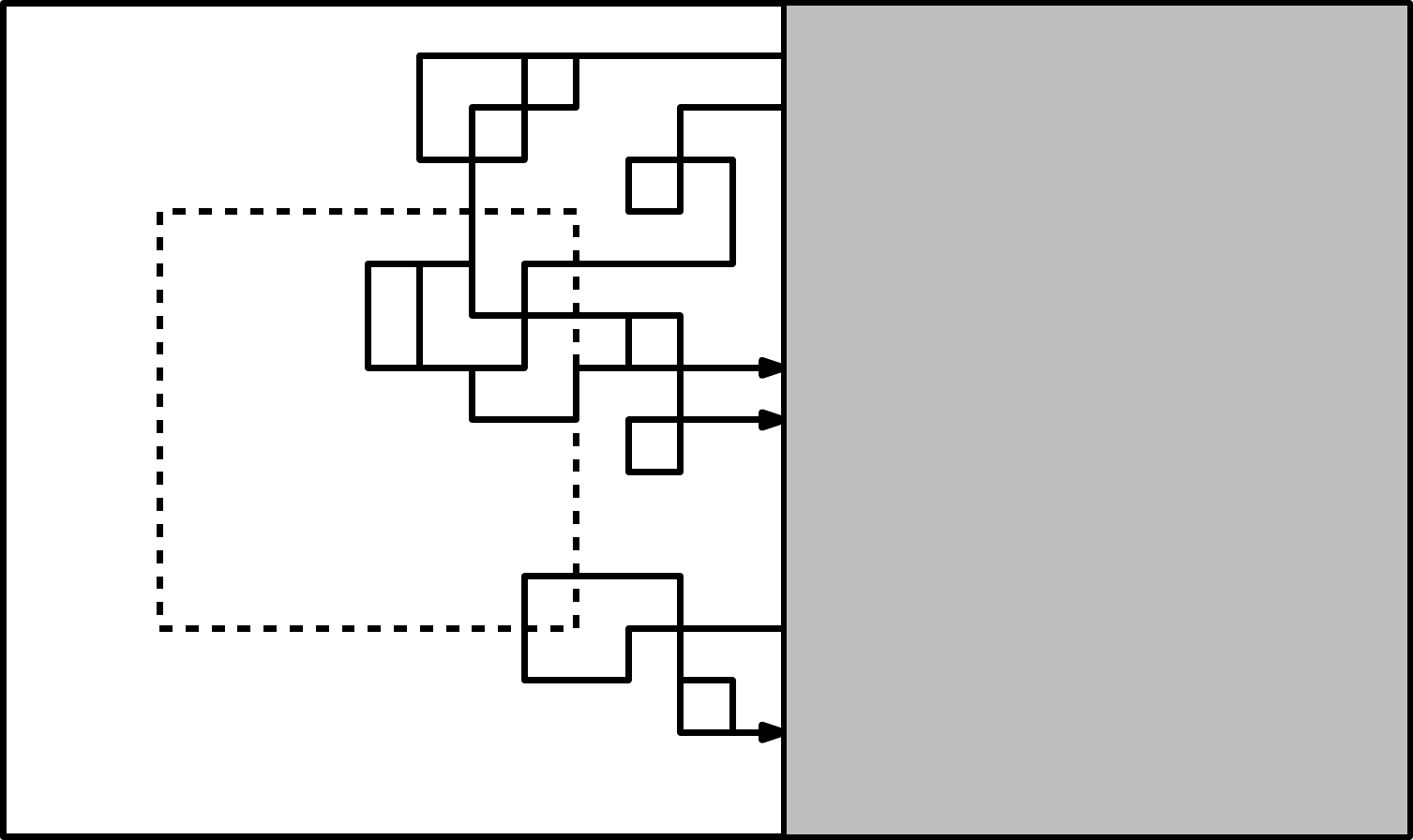}
\hskip 5mm
\includegraphics[scale=0.5]{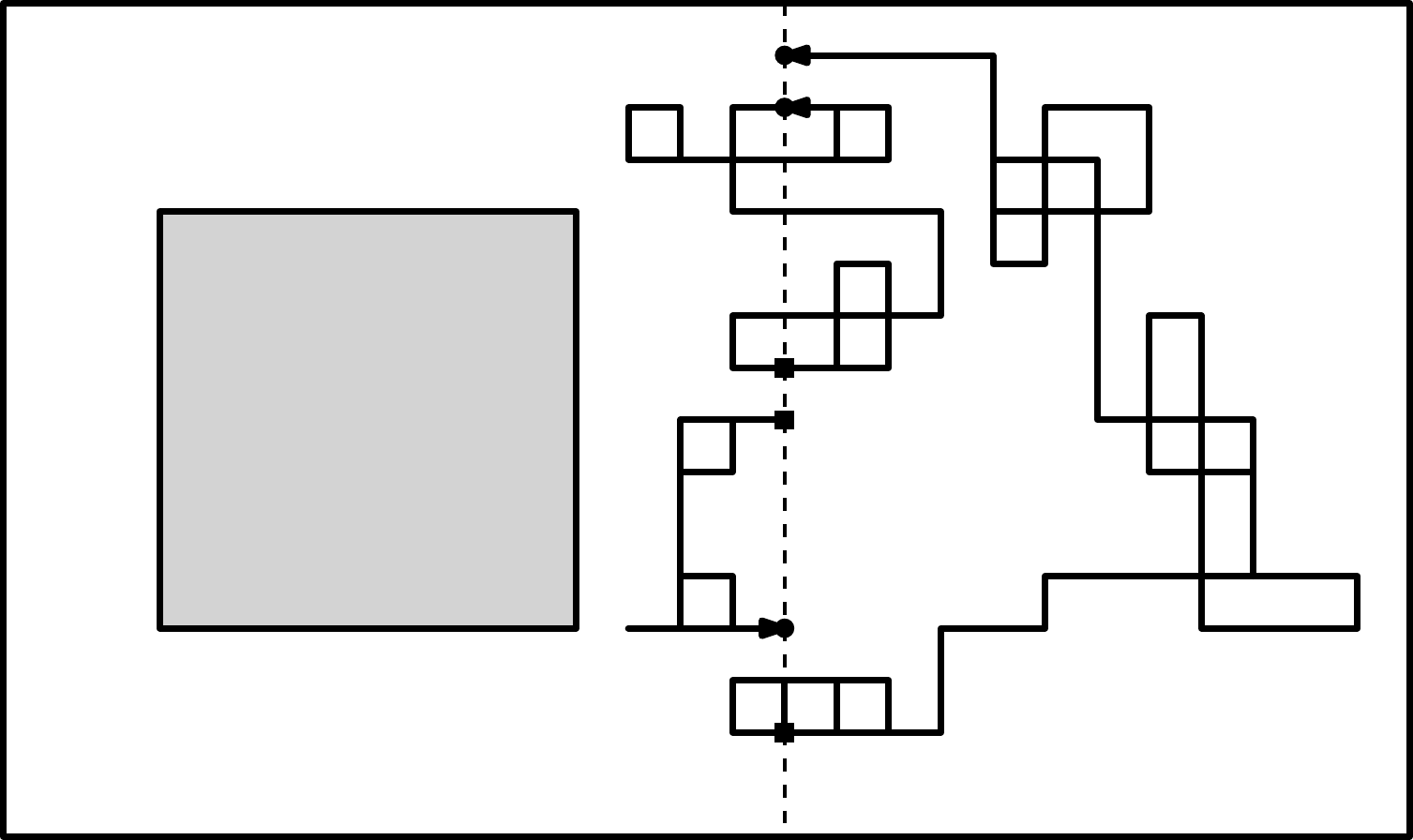}
\end{center}
\caption{Discovering (i) the oriented excursions away from the right part that reach the small square, (ii) sampling the three oriented bridges in the complement of the small square.}
\end{figure}

Note that this conditional distribution is fully described by the vectors $\X$ and $\Y$ (ie. it depends on $\eta$ just as a function of $\X$ and $\Y$), which is one of the main features of this result. In other words, conditionally on $\X$ and $\Y$, $\eta$ and $\beta$ are independent. In particular, the number of actual loops that are being created by $\beta$ when one concatenates it with $\eta$ does not intervene in the conditional distribution, which is a specific feature of this $\alpha=1$ case.

Let us comment on the case where $F_2 = D \setminus F_1$: If one then conditions on the number of jumps of the loop-soup on each edges from a point in $F_1$ to a point of $F_2$ (one then gets a collection $(\X_j', \X_j)_{j \le \N}$ of jumps from $\X_j' \in F_1$ to $\X_j \in F_2$), and on the number of jumps of the loop-soup on each edge from a point of $F_2$ to a point of $F_1$ (one then gets a collection $(\Y_j, \Y_j')_{j \le \N}$ of jumps from $\Y_j \in F_2$ to 
$\Y_j' \in F_1$), then the conditional distribution of the 
 missing pieces in $F_2$ and in $F_1$ are independent, and there are 
respectively the unordered bridge measure in $F_2$ from $\X$ to $\Y$ (this corresponds to $\beta$), and the unordered bridge measure from $\Y'$ to $\X'$ in $F_1$ (this corresponds to $\eta$ without the first and last jumps of each excursion). This can be interpreted as a spatial Markov property of the occupation field on oriented edges (the random function that assigns to each oriented edge the total number of jumps of the soup along this edge) of the $\alpha=1$ soup of oriented loops. We will discuss this again at the end of this section. 

In the same spirit, we can in fact ``symmetrize'' also Proposition \ref {p1} also when $F_2$ is a subset of the complement of $F_1$. Let us then define the collection of crossings $\eta_{1 \to 2}$ to be the parts of the loops in the loop-soup of the type $a_0, e_1, \ldots, a_n$ with $a_0 \in F_1$, $a_n \in F_2$ and $a_1 , \ldots, a_{n-1} \in D\setminus (F_1 \cup F_2)$. 
We also define $\eta_{2 \to 1}$ similarly, and note that there are as many crossings from $F_1$ to $F_2$ as there are crossings from $F_2$ to $F_1$. Let 
$\X$ (resp. $\X'$) denote the vector of endpoints of $\eta_{1 \to 2}$ (resp. $\eta_{2 \to 1}$) and $\Y$ (resp. $\Y'$) the vector of starting points of $\eta_{2 \to 1}$ (resp. $\eta_{1 \to 2}$). Then, we can note that $\X$ and $\Y$ are exactly the same as the ones defined in Proposition \ref {p1}, while $\X'$ and $\Y'$ correspond to those that one obtains when interchanging $F_1$ and $F_2$. Furthermore, $\eta_{1 \to 2}$ and $\eta_{2 \to 1}$ are fully determined by $\eta$ (or alternatively by the symmetric family $\eta'$
of excursions outside of $F_1$ that do reach $F_2$). It follows readily from Proposition \ref {p1} that: 

\begin {proposition}
\label {p1bis}
 Conditionally on $\eta_{1 \to 2}$ and on $\eta_{2 \to 1}$, the missing parts of the loops that they are part of (these are the loops of the $\alpha=1$ soup of oriented loops that intersect both $F_1$ and $F_2$) are described by two independent unordered bridges with conditional distributions $B_{\X,\Y}^{D \setminus F_1}$ and $B_{\X', \Y'}^{D \setminus F_2}$. 
\end {proposition}

Note that the other loops in the loop-soup (i.e. the loops that either do not intersect at least one of the two sets $F_1$ or $F_2$) are just described by a loop-soup in the complement of $F_1$ and a loop-soup in the complement of $F_2$, that are coupled to share exactly the same loops that stay in $D \setminus (F_1 \cup F_2)$. 

\medbreak

Let us now prove Proposition \ref {p1}. 
\begin {proof}
Let us consider a family $E$ of $N$ excursions, such that $P (\eta = E) >0$ and such that the $N$ excursions $E_1, \ldots , E_N$ of $E$ are all different. Then if $\eta=E$ and $\tilde \LL = \overline {L}$, all the loops in $\overline L$ are simple, and they do occur necessarily exactly once (and not more). Hence, for such an $\overline L$, the probability that  $\overline {\mathcal L} = \overline L$ is proportional to $g^{-\overline n (\overline L) }$ where $\overline n$ is the sum of the lengths of the loops in ${\overline L}$ (and the proportionality constant does not depend on $\overline L$). 

 On the other hand, if $X$ and $Y$ are the vector of end-points of $E$, the $B_{X,Y}^{D \setminus F_1}$-probability to sample a unordered bridge that gives rise exactly to $\overline L$ when concatenating it to $E$ is proportional  to  $g^{- K}$ (where $K = \overline n (\overline L) - \overline n (E)$ is the total length of the generalized bridge), because there is just one permutation per bridge that works.  It therefore follows immediately that conditionally on ${\eta=E}$, the distribution of the missing bridges is indeed $B_{\X,\Y}$ in $D \setminus F_1$.
 
 Instead of treating directly the case of multiple occurrences of the same excursions in $\eta$, we will use the following trick (a similar idea can be used to show the fact that the loops erased during Wilson algorithm do correspond exactly to an oriented loop-soup, see for instance \cite {W1}). We choose a very large integer $W$ (that is going to tend to infinity), and we decide to replace the graph $\G$ by the graph $\G^W$, which is obtained by keeping the same set of vertices as $\G$, but where each edge of $\G$ is replaced by $W$ copies of itself. In this way, each site has now $gW$ outgoing edges instead of $g$. There is of course a straightforward relation between random walks, loops and bridges on $\G^W$ and on $\G$. For instance, a loop-soup (resp. bridge, resp. excursion) on $\G^W$ is directly projected on  a loop-soup (resp. bridge, resp. excursion) on $\G$. 

Let us couple loop-soups with intensity $\alpha=1$ in all of the $\G^W$'s on the same probability space, in such a way that the projections of the loop-soups in $\G^W$ onto $\G$ (in the sense described above) are the same for all $W$'s. We fix also $F_1$, $F_2$, and define (with obvious notation), $ \LL^W$, $ \eta^W$, $\overline L^W$ etc. Note that the vectors of extremal points $\X$ and $\Y$ are then the same for all $\eta^W$'s.

We can also note that the probability that some edge is used more than once in the loop-soup does tend to $0$ as $W \to \infty$.  
The probability that all excursions in $\eta^W$ are different therefore tends to $1$ as $W \to \infty$. 

But conditionally on the fact that all excursions in $\eta^W$ are different (applying our previous result to $\G^W$), we know that the conditional distribution of $\overline \LL^W \setminus \eta^W$ given $\eta^W$ is the 
bridge probability measure from $\X$ to $\Y$ in $D^W \setminus F_1$. Projecting this onto $\G$, we get that the conditional distribution of $\beta$ given $\eta^W$ (on the event that in $\eta^W$, no two excursions are the same) is the unordered bridge measure $B_{\X, \Y}$ in $D \setminus F_1$. 

If $U(W)$ is the event that no two excursions of $\eta^W$ appear twice, we therefore get that, conditionally on $\eta = E$ and $U(W)$, the conditional distribution of $\beta$ is the unordered bridge measure $B_{\X,\Y}$ in $D \setminus F_1$. We now just let $W \to \infty$, which concludes the proof of the proposition. 
\end {proof}

\subsection {Partial resampling of  soups of unoriented loops at $c=1$}

Let us now come back to the setting where the graph $\G$ is unoriented. 
When one considers a soup of unoriented loops with intensity $\nu$ (recall that this corresponds to $c=1$ or $\alpha = 1/2$ ie. to a soup of oriented loops with intensity $\mu/2$ where we forget the orientation of each loop). We denote the collection of unoriented loops that intersect both $F_1$ and $F_2$ by 
${\overline \LL}= ( \tilde \LL_1, \ldots, \tilde \LL_\M)$, the corresponding collection of (unoriented) excursions by $\eta=  ( \tilde \eta_1, \ldots, \tilde \eta_\N)$   and the endpoints of these $\N$ excursions by $\Z= (\z_1, \ldots, \z_{2\n})$. The missing parts of the (unoriented) loops are unoriented paths that join each $\z_i$ to exactly one other $\z_j$, so that $\beta$ is an unordered $\Z$-bridge in $D \setminus F_1$. 

Note again that it is intuitively possible to explore the excursions $\tilde \eta_j$ ``starting'' from their intersections with $F_1$ in both directions, until hitting $F_2$ (and in this way, one did yet discover 
the missing parts $\beta$).

\begin {proposition}
\label {p2}
 The conditional distribution of $\beta$ given $\eta$ is exactly the unordered unoriented bridge measure $B_{\Z}$ in $D \setminus F_1$.
\end {proposition}

\begin{figure}[ht]
\begin{center}
\includegraphics[scale=0.5]{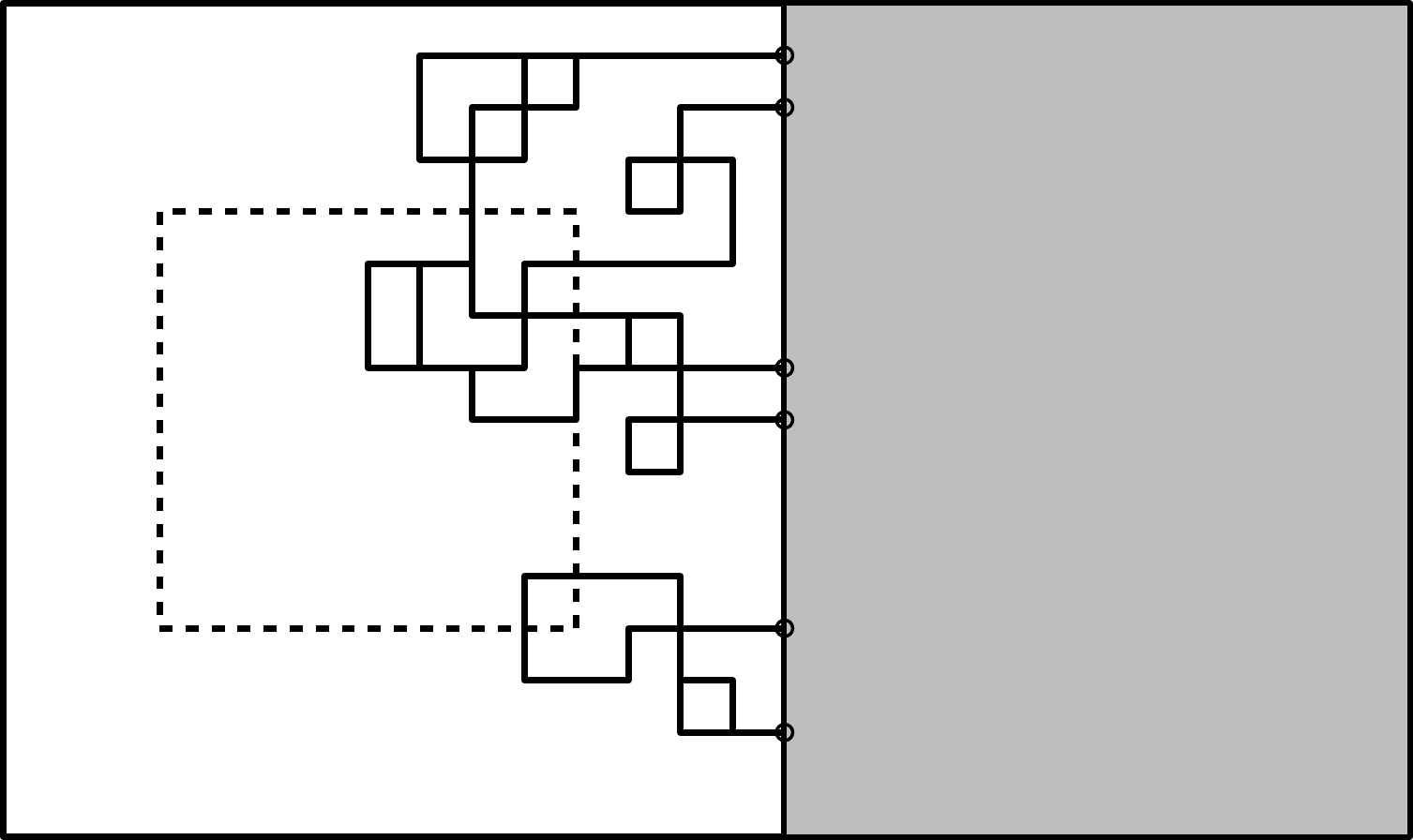}
\hskip 5mm
\includegraphics[scale=0.5]{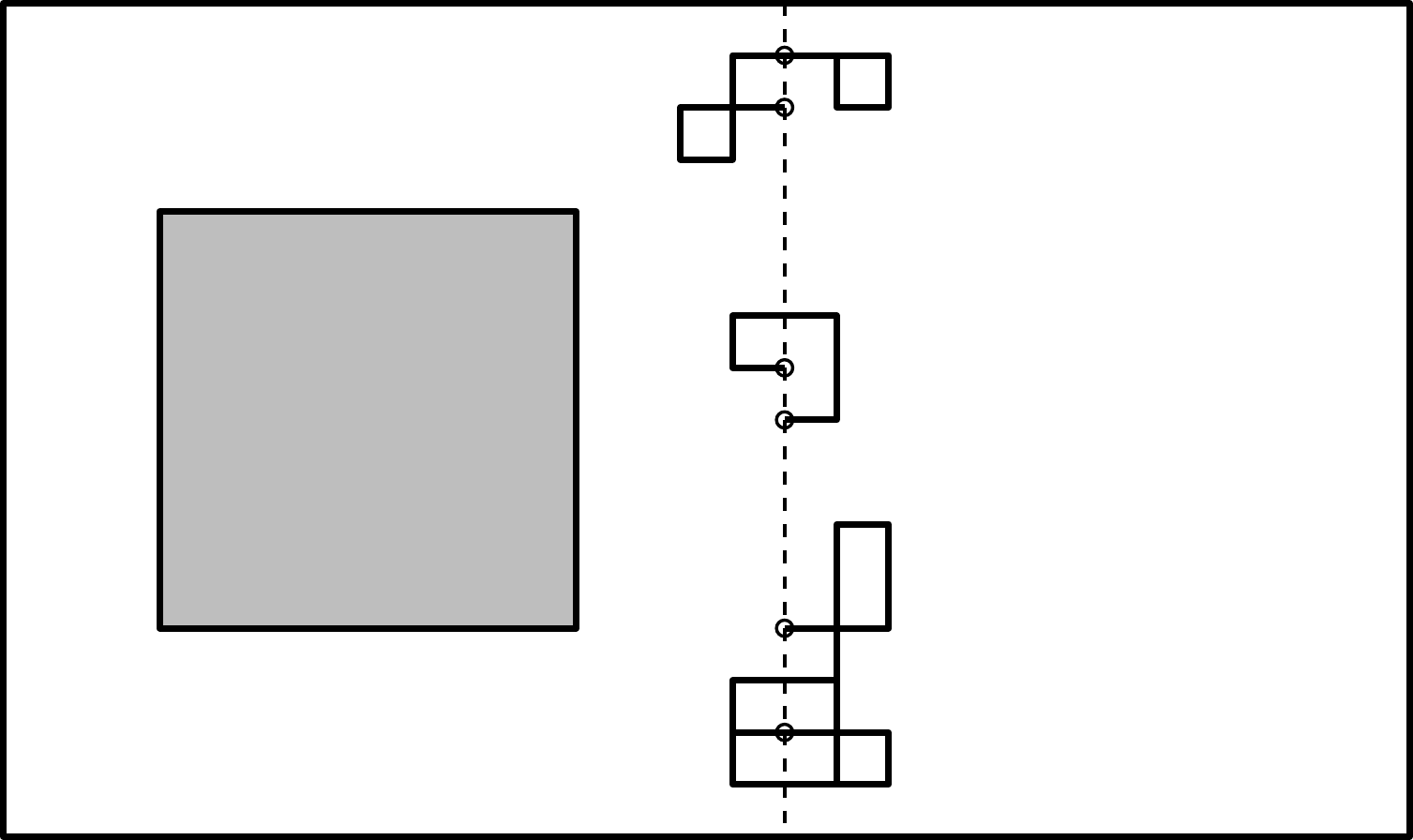}
\end{center}
\caption{Discovering (i) the unoriented excursions away from the right part that reach the small square,  (ii) sampling the three unoriented bridges in the complement of the small square.}
\end{figure}

Just as in the oriented case, we stress that an important feature in this statement is that this conditional distribution is a measurable function of the vector $\Z$ (the other information on the excursions are not needed).
 We will further comment on this in the next subsection.

\begin {proof}
We will follow the same idea as in the proof of the oriented case.
 As in the unoriented case, when the $N$ pieces $\tilde E_1, \ldots , \tilde E_N$ of $E$ are all different, the statement is almost immediate (for each good ordered bridge, only one pairing works in order to complete $E$ into $\overline L$, and the probability to complete these $N$ pieces into $\overline L$ is therefore proportional to $g^{-K}$ where $K$ is the difference between the total number of jumps in the loop-configuration and in $E$).
 
We then use the same trick with copying each edge a large number of times. The very same argument the works, almost word for word.
\end {proof}

\subsection {Spatial Markov properties} 

The particular case  where $F_2$ is the complement of $F_1$ is also of interest for the soup of unoriented loops. Let us for instance describe how things work for the occupation times of loop-soups (which is the main focus of the papers of Le Jan \cite {LJ}). If one then conditions on the numbers of jumps of the loop-soup on all edges between a point in $F_1$ and a point of $F_2$ (in either direction -- the loops being unoriented there is anyway no direction), then the conditional distribution of the parts $\beta$ in $F_2$ of the loops that intersect both $F_1$ and $F_2$ is described by Proposition \ref {p2} and it is a unordered unoriented bridge in $F_2$ (and it is in fact fully described by the knowledge of the number of jumps along the edges between $F_1$ and $F_2$, ie. this conditional distribution is a function of these number jumps of the edges between $F_1$ and $F_2$). But, the situation is symmetric and we can interchange the roles of $F_1$ and $F_2$; we therefore conclude that given $\beta$ and the 
numbers of jumps along the edges between $F_1$ and $F_2$, the conditional distribution of $\eta'$ defined to be the collection $\eta$ where one has removed the two extremal jumps of each $\eta_j$ (these are the jumps between $F_1$ and $F_2$), is that of an unordered unoriented bridge in $F_1$ (and the law of this bridge is also fully described by the number of jumps between $F_1$ and $F_2$). 
    
In other words, when one conditions on these number of jumps along the edges  between $F_1$ and $F_2$, we can enumerate these jumps (using some deterministic lexicographic rule) by $(\Z_j', \Z_j)_{j \le 2 \N}$ where $\Z_j' \in F_1$ and $\Z_j \in F_2$. Then, the conditional distribution of $\eta'$ and $\beta$ are conditionally independent unordered bridges, respectively following the unordered bridge measures $B_{\Z'}^{D \setminus F_2}$ and $B_{\Z}^{D \setminus F_1}$. In particular, when adding on top of this the loop-soups in $F_1$ and the loop-soups in $F_2$, it follows that conditionally on the occupation times (i.e. on the number of jumps $N_e$ across each edge) on the edges between $F_1$ and $F_2$, the occupation times on sites and edges in $F_1$ is independent of the occupation times on sites and edges in $F_2$. We can rephrase this property in the following sentence: 
The occupation time field on edges of the soup of unoriented loops for $c=1$ does satisfy the spatial Markov property.

We can note that if $U$ is a non-negative function of the occupation time field on the edges of the form $U((N_e)) = \prod_e u_e ( N_e)$, such that the expectation of $U$ (for the $c=1$ loop-soup) is equal to one, then if we define the new probability measure $Q$ on occupation times on edges by $dQ / dP ((N_e)) = U((N_e))$, then the spatial Markov property also holds for $Q$. This can be used to represent a modification of the Markov chain (ie. different walks with non-uniform jump probabilities). 

If we consider an unoriented graph, but that we interpret as an oriented graph (each unoriented edge defines an oriented edge in each direction), on which we define an $\alpha=1$ soup of oriented loops, then we can also reformulate the results of Subsection \ref {ss2.1} in a similar way. More precisely, for each edge, we can define the total  number of jumps $N_1 (e)$ by the soup in one direction of $e$, and $N_2 (e)$, the number of jumps in the opposite direction. Then, if we define  $N_e := ((N_1 (e), N_2 (e))$, this 
two-component occupation time field on edges of the $\alpha=1$ soup of oriented loops satisfies the spatial Markov property in the same sense as above.

\medbreak
Let us now come back to the study of the loops themselves, and not just of the cumulated occupation time of the soup. 
As in the oriented case, we can also (when $F_2$ is a subset of $F_1$) rephrase Proposition \ref {p2} in a more symmetric way, involving the crossings between $F_1$ and $F_2$. We define 
$\eta_{1 \leftrightarrow 2}$ the set of (unoriented) parts of loops  in the $c=1$ loop-soup that join a point of $F_1$ to a point of $F_2$ and otherwise stay in the complement of $F_1 \cup F_2$, and we denote by $\Z$ the vector of endpoints of these crossings in $F_2$, and by $\Z'$ the set of endpoints in $F_1$. Then: 

\begin {proposition}
\label {p3bis}
Conditionally on $\eta_{1 \leftrightarrow 2}$, the missing parts of the unoriented loops that these crossings are part of (these are the loops in the loop-soup that intersect both $F_1$ and $F_2$) are described by two independent unordered unoriented bridges with respective conditional distributions $B_\Z^{D \setminus F_1}$ and $B_{\Z'}^{D \setminus F_2}$.
\end {proposition}
Figure \ref {Bexc} that illustrates the corresponding result in the Brownian case, can also be used to illustrate this result. 

It is also easy to generalize Proposition \ref {p3bis} and Proposition \ref {p1bis} to more than two sets $F_1$ and $F_2$ (and have instead $n$ disjoint sets $F_1, \ldots , F_n$). For instance, in the unoriented case, one then conditions on the set $\eta_\leftrightarrow$ of all crossings from 
any $F_i$ to any other $F_j$ that also stay in the complement of all the other $F_k$'s. These crossings define $n$ vectors $\Z^1, \ldots, \Z^n$ (where $\Z^j$ is a list of the even number of endpoints on $F_j$ of the aforementioned crossings). Conditionally on $\eta_\leftrightarrow$, the missing parts of the loops (that are the loops in the loop-soup that touch at least two different $F_j$'s) are described by $n$ conditionally independent unordered unoriented bridges with respective distributions $B_{\Z_j}^{D' \cup F_j}$ (where $D = D \setminus \cup_i F_i$) for $ j \le n$.

\medbreak

Such decompositions of the loops in the soup that intersect  disjoint compact sets into crossings +  conditionally independent unordered bridges, can be immediately transcribed to the case of Brownian loops on the cable system associated to this graph as studied in \cite {Lupu1}; we leave this as a simple exercise to the reader.
This is all of course closely related to the Markov property of the Gaussian Free Field, as well as to Dynkin's isomorphism theorem \cite {Dy} via the relation between the square of the GFF and the loop-soup (see e.g.. \cite {LJ} and the references therein for background).  

\medbreak
With such Markovian-type properties in hand, a natural next step is to define random sets that play the role of stopping times for one-dimensional Markov processes. In the setting of the discrete GFF, these are the local sets as defined in \cite {SS}, and that turned out to be very useful concepts. Just as for one-dimensional stopping times, there are several possible ways to define them, depending on what precise filtration on considers. 
In the present case (we do here describe the definitions in the unoriented loop-soup for $c=1$, but the oriented case would be almost identical), one can for instance say that: 
\begin {itemize}
 \item A random set of points ${\mathcal  F}$ is a stopping set for the occupation time field filtration, if for any $F_1$, the event $\{ {\mathcal  F} =F_1\}$ is measurable with respect to the occupation time field on all edges adjacent to $F_1$.
 \item A random set of points ${\mathcal  F}$ is a stopping set for the loop-soup filtration, if for any  $F_1$, the event $\{ {\mathcal  F} =F_1\}$ is measurable with respect to the trace of the loop-soup on all edges adjacent to $F_1$ (i.e. it is measurable with respect to the set of loops that are fully contained in $F_1$ and the set of excursions $\eta$ defined above, when $F_2$ is the complement of $F_1$.
 \item  A random set of points ${\mathcal  F}$ is a stopping set for the loop-soup, if for any  $F_1$ such that $P ({\mathcal  F}=F_1 ) >0$, conditionally on the event $\{{\mathcal  F} =F_1\}$, the distribution of the loop-soup 
 outside of $F_1$ consists of the union of an independent loop-soup in the complement $F_2$ of $F_1$ and of a set of bridges in $F_2$, with law described as above via the end-points of the excursions $\eta$ in $F_1$. 
\end {itemize}
Clearly, the first definition implies the second one, which implies the third one by Proposition \ref {p2} (the third property for the first two definitions can be viewed as a ``strong Markov property'' of these fields), but the converse is not true (the last definition allows the use of ``external randomness'' in the definition of ${\mathcal  F}$ (while the second does not), and the second one allows features of individual loop (while the first does not).

\subsection {Brownian loop-soup decompositions}

The previous results have almost identical counterparts in the setting of oriented Brownian loop-soups with intensity $\alpha=1$ and unoriented Brownian loop-soups with intensity $c=1$. 
 
Suppose that $D$ is an open subset of $\R^d$ such that the (Dirichlet) Green's function in $D$ is finite (away from the diagonal). 
Suppose that $F_1$ and $F_2$ are two disjoint compact sets in $D$, that are both non-polar for Brownian motion (i.e. Brownian motion started away from these sets has a non-zero probability to hit them). 
Then, we can again define: 
\begin {enumerate}
 \item 
 The law of unordered oriented Brownian bridges in $D \setminus F_1$ from a finite family $X=(x_1, \ldots, x_n)$ of points to another such family $Y=(y_1, \ldots, y_n)$,
and the law of unordered unoriented $Z$-Brownian bridges in $D \setminus F_1$ from a finite family of points
$Z = (z_1, \ldots , z_{2n})$ to itself (in the latter case, points of $Z$ are paired, like in the random walk case). This works as long as all Green's functions involved are finite (which is the case as soon as all $x_i\not=y_j$ for all $i,j$, and that $z_i \not= z_j$ for all $i \not=j$). 
\item
 The set $\eta$ of $\N$ oriented (resp. unoriented) excursions of the loops in an oriented (resp. unoriented) loop-soup with intensity $\alpha=1$ (resp. $c=1$) away from $F_2$, that reach $F_1$. In the ordered case, we 
call their endpoints vector $\X = (\X_1, \ldots, \X_\N)$ and their starting point vector $\Y$,  and in the unoriented case, we call $\Z= ( \Z_1, \ldots, \Z_{2\N})$ the extremity vector. 
\end {enumerate}
Then, the Brownian counterparts of Proposition \ref {p1} and of Proposition \ref {p2} go as follows: 
\begin {proposition}
- For the soup of oriented Brownian loops with $\alpha=1$: Conditionally on $\eta$, the missing pieces of the loops (that the pieces $\eta$ are part of) are distributed like an unordered Brownian bridge from $\X$ to $\Y$ in $D \setminus F_1$.

- For the soup of unoriented Brownian loops with $c=1$: Conditionally on $\eta$,
the missing pieces of the loops are distributed like an unordered unoriented $\Z$-Brownian bridge in $D \setminus F_1$.
\end {proposition}

\begin{figure}[ht!]
\begin{center}
\includegraphics[scale=0.55]{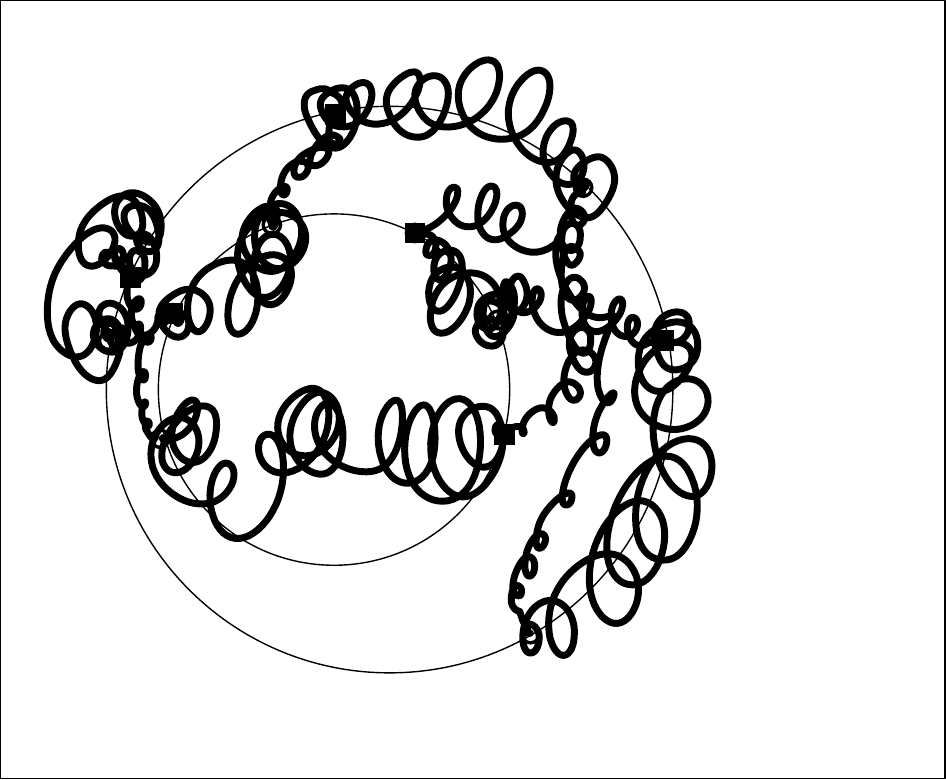}
\quad
\includegraphics[scale=0.55]{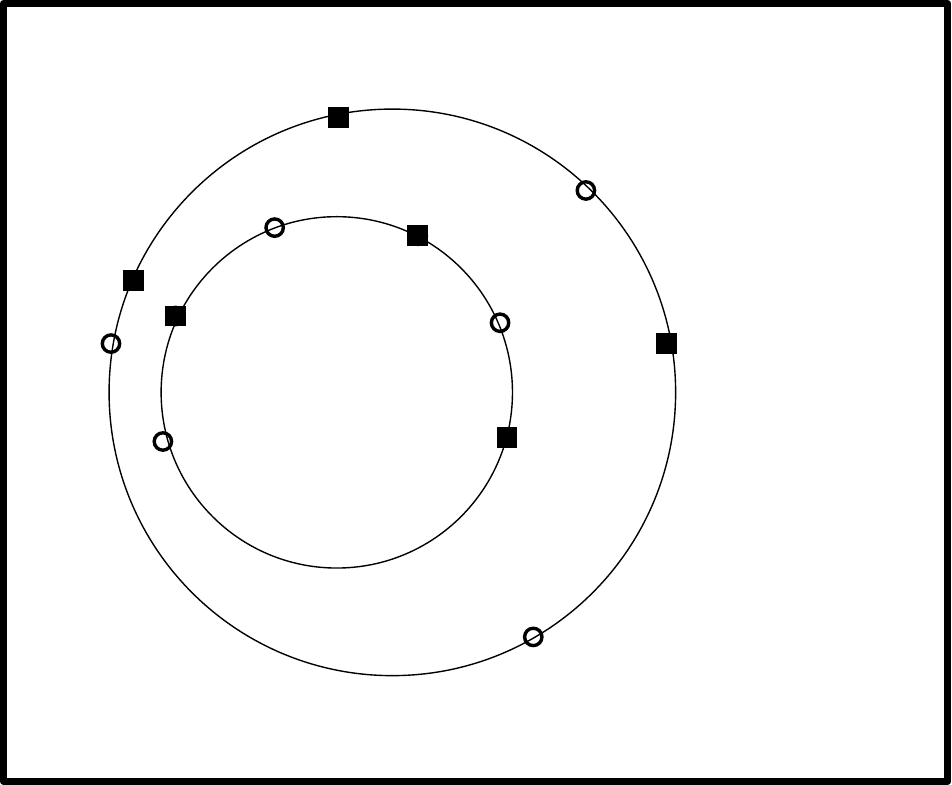}
\end{center}
\caption{Sketch of the oriented Brownian case: (i) The two oriented loops that touch the two circles, (ii) keeping only the endpoints of these crossing on each circle, with trace of the orientation}
\label {Bexc}
\end {figure}

\begin {figure}[ht!]
\begin {center}
\includegraphics[scale=0.55]{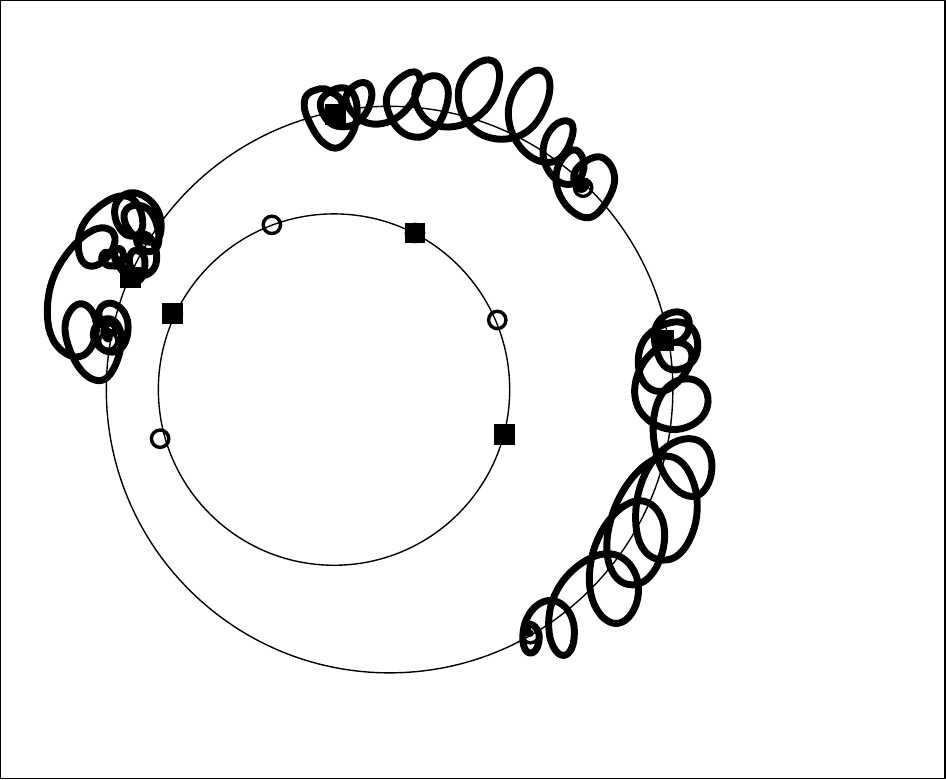}
\includegraphics[scale=0.55]{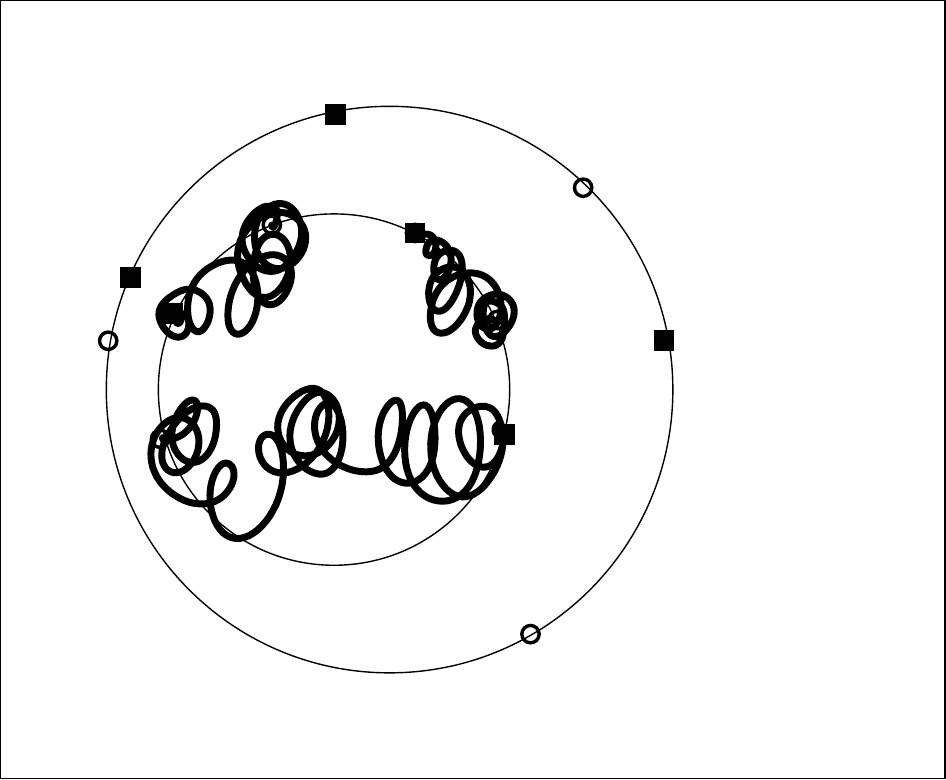}
\includegraphics[scale=0.55]{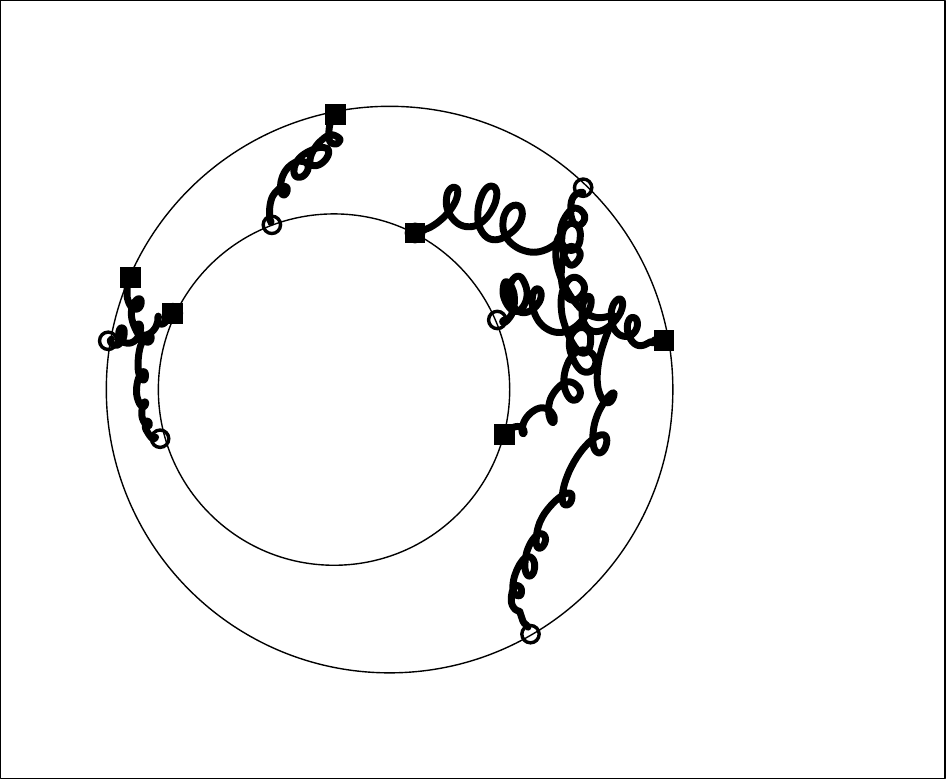}
\end {center}
\caption{(iii) The outer bridges joining each circle point to a square point, (iv) sampling the inner bridges joining each circle point to a square point, (v) the six crossings, joining a circle point to a square point. The final loops are oriented so that the crossings from small to large circle go from a circle point to a square point}
\label {Bexc2}
\end{figure}

And as before, one can derive the more symmetric results:  
For instance, if $F_1$ and $F_2$ are two disjoint compact subsets of $D$, we can define the crossings from $F_1$ to $F_2$ and vice-versa in the oriented case, and the crossings between 
$F_1$ and $F_2$ in the unoriented case. When one conditions on these crossings, one can then complete the picture with two conditionally independent unordered oriented bridges (in the oriented case) or by two conditionally independent unordered unoriented bridges (in the unoriented case). We illustrate this  result in Figures \ref {Bexc} and \ref {Bexc2} (here we consider the oriented case, $D$ is the rectangle, $F_1$ is the small circle and $F_2$ the large circle). Conditionally on the points (and their status -- square or circle depending on the orientation of the loops) on the two circles, the three pictures in Figure \ref {Bexc2} are independent (this is the oriented version of Figure \ref {f2}).  

\medbreak

In the context of two-dimensional continuous systems, clusters of loops in a loop-soup are interesting to study, as  pointed out in \cite {Wcras};
it has been proved in \cite {ShW} that boundaries of such clusters for $c \le 1$ form Conformal Loop Ensembles with parameter $\kappa = \kappa(c)$, where $\kappa (1) =4$. The CLE$_4$ (and the SLE$_4$ curves) is also known (see \cite {SS,Du}) to be related quite directly to the Gaussian Free Field. 
The role of the $c=1$-clusters of loops in the framework of cable-systems and in relation to the Gaussian Free Field has been pointed out by Lupu \cite {Lupu1} (the clusters  provides a direct link between the loop-soups and the Gaussian Free Field itself, rather than just to its square).
The present result sheds some light on the recently derived \cite {QW} decomposition of critical $2d$ loop-soup clusters (for $c=1$) in terms of Poisson point processes of Brownian excursions (we refer to \cite {QW} for comments and questions).

\section {Resampling for continuous-time loop-soups, the GFF and random currents}

We devote now a short separate section on the case of discrete continuous-time loop-soups, that have been studied by Le Jan \cite {LJ}. As we shall see, in that setting, it is 
natural to consider the conditioned distribution of the loop-soup (unoriented for $c=1$ ie. $\alpha = 1/2$, or oriented for $\alpha=1$) given the value of their local times on
a given family of sites. Some of the results are very closely related to Dynkin's isomorphism theorem (ie. it will be a pathwise version of a generalization of it). 
Just as previously, we will describe the case of simple random walk on the graph where each point $x$ has the same number $g$ of outgoing edges, but the results can easily be generalized to the case of general Markov chains. Some of following considerations will be reminiscent of the arguments in \cite {LJ} (sections 7 and 9 in particular). 
In the first subsections, we will focus on the case of unoriented loop-soups, and we will briefly indicate the similar type of results that one gets in the oriented case. 

\subsection {Slight reformulation of the resampling property of the discrete loop-soup} 

We can start with the same setting as before, with the graphs $\Gamma$ and $D \subset \Gamma$, the random walk on this graph killed upon hitting $\Gamma \setminus D$, and its Green's function $G_D ( \cdot, \cdot)$. 
In the previous sections, we chose for expository reasons (as this was for instance the natural preparation for the Brownian case) to study loops in the loop-soup that visit two different sets of sites $F_1$ and $F_2$. 
But in fact, the following setting is a little more natural and more general: Consider now a family $e_1 , \ldots, e_n$ of edges of $D$, and the graph $D'$ obtained by removing these $n$ edges from $D$.
We can now sample an unoriented loop-soup (for $c=1$), and observe the numbers $N_1, \ldots, N_n$ of jumps along thoses $n$ unoriented edges. We now want to know 
the conditional distribution of the entire loop-soup given this information. In particular, we would like to know how these $N_1 + \ldots + N_n$ jumps are hooked together into loops (clearly, the loop-soup in $D'$ consisting of the loops that use none of these $n$ edges is independent of $ \overline N := (N_1, \ldots, N_n)$). 

We can associate to $ \overline N$ the vector $\Z$ consisting of the $2N_1 + \ldots + 2 N_n$ endpoints of these jumps. Once we label them, we can as before the collection $\beta$ of pairing and bridges that join them in the loop-soup. Note that the bridge is allowed to contain no jump when one pairs two identical end-points. We can also define the unordered bridge measures in $D'$ (corresponding to paths that use no edge of $D \setminus D'$) as before. Then, exactly as before, one can prove the following version of the resampling: 
\begin {proposition}
\label {p5}
 The conditional distribution of $\beta$ given $N_1, \ldots , N_n$ is exactly the unordered unoriented bridge measure $B_{\Z}$ in $D'$.
\end {proposition}
Note that for some choices of family of edges $e_1, \ldots, e_N$, it can happen that an even number of endpoints of the discovered jumps are at a certain vertex where no neighboring edge is in $D'$. In that case, the bridge measure pairs these jumps at random and the corresponding bridge is anyway the empty bridge from $x$ to $x$. A trivial example is of course the case where $e_1, \ldots, e_n$ are all the edges of $D$. Then, the proposition just says that the conditional distribution of the loops given the occupation time measure is obtained by just pairing at random the incoming edges at each site. ``Loops can exchange their hats uniformly at random at each site''. 

This reformulation makes it clear that in the discrete time setting, the Markov property of the occupation time field is really a Markov property on the edges (which is not surprising, given that the field is actually naturally defined on the edges). 

\subsection {Continuous-time loops}

Following  Le Jan's approach  \cite {LJ}, we now introduce the associated continuous-time Markov chain, for each site $x$, the chain stays an exponential waiting time of mean $1/g$ before jumping along one of the $g$ outgoing edges chosen at random (for expository reasons, we describe this in the case where each edge has the same number of outgoing edges). Note that we allowed stationary edges in the graph, so that 
the continuous-time Markov chain can also ``jump'' along those (and we can keep track of these jumps, even if they do not affect the occupation time at sites). As pointed out by Le Jan, the loop-soup of such continuous-time loops for $\alpha=1/2$ is particularly interesting, as its cumulated occupation time (on sites) is exactly the square of a Gaussian Free Field on this graph (here one may introduce one or more killing point, so that the loop-soup occupation-time is finite, and the free field with boundary value $0$ at this point is well-defined). 
In this setting, the loops of the discrete Markov chain do correspond exactly to loops of the continuous-time chain, but the latter also contains some additional stationary loops,
that just stay at one single point without jumping during their entire life-time. 

When one considers a continuous-time loop and a finite set of vertices in the graph that it does visit, one can cut-out from the loops the time that it does spend at these points and obtain a finite sequence of excursions away from this set. This corresponds 
to the usual excursion theory of continuous-time Markov processes (an excursion from $x$ to $y$ will be a path that jumps out of $x$ at time $0$ and jumps into $y$ at the endpoint of the excursion). One can the introduce the natural excursion measure 
$\mu_{x,y}^A$, which is the natural measure on set of unoriented
excursions that go from $x$ to $y$ while avoiding all the points in $A$ (it corresponds to the discrete excursion measure that puts a mass $g^{-n}$ to such an excursion with $n$ jumps, and one then adds $n-1$ independent exponential waiting times at the $(n-1)$ points inside the excursions.  

One can 
view the continuous-time Markov chain as the limit when $M \to \infty$ of the discrete-time Markov chain on a graph $D^M$, where one has added to each site $x$,  $M$ stationary edges from $x$ to itself (when one renormalizes time by $1/M$, the geometric number of successive jumps along these added stationary edges from $x$ to $x$ before jumping on another edge, does converges to the exponential random variables) -- this approach is for instance used in \cite {W1} in order to derive the properties of the continuous-time chains and loop-soups from the properties of the discrete-time loop-soups. 
Let us now consider a finite set of points $x_1 , \ldots, x_n$ in the graph, and for a given $M$, we condition on the  $N_1, \ldots, N_n$ of jumps by the loop-soup along the stationary unoriented edges $e_1, \ldots, e_n$. More precisely $N_1$ will denote the total number of jumps in the loop-soup along the $M$ added stationary edges from $x$ to $x$. Note that because both end-points of a stationary edge are the same, these $N_1$ jumps correspond to $2N_1$ jump-endpoints, that are all at $x_1$.  
We can now apply Proposition \ref {p5} to  this case; this describes the distribution of how to complete and hook up these $N_1 + \ldots + N_n$ jumps into unoriented loops in order to recover the loops in the loop-soup that they correspond to.
One has to pair all these $2N_1 + \ldots + 2N_n$ endpoints.

Mind that as $M$ gets large, the mass of the trivial excursion from $x_1$ to $x_1$ with zero life-time is always $1$, while the mass of (unoriented) excursions with at least one jump along the
``non-added'' $nM$ stationary edges neighboring these points 
from $x_1$ to some $x_j$ that stays away from $\{x_1, \ldots, x_n \}$ during the entire positive lifetime (if it is positive) will be of order $1/M$ (unless all neighbors of $x_1$ are in $\{x_1, \ldots, x_n \}$ in which case this 
quantity is zero) and 
that the set of excursions from $x_1$ to $x_j$ that visit at least one of the points of $\{x_1, \ldots, x_n \}$ during its positive life-time is of the order of $O(1/M)^2$. 
It is a simple exercise that we safely leave to the reader to check that in the $M \to \infty$ limit, the discrete Markovian description becomes the following: 

\begin {proposition}
 If we consider the continuous-time Markov chain loop-soup and condition on the total occupation time $l(x_1), \ldots, l(x_n)$ at the $n$ points $x_1, \ldots, x_n$, then the unoriented 
 excursions away from this set of points 
 by the loop-soup will be distributed exactly like a Poisson point process of excursions with intensity $\mu_l = (1/2) \times \sum_{ i \le j}   l(x_i) l(x_j) \mu_{x_i, x_j}^{x_1, \ldots, x_n} $
 conditionned on the event that the number of excursions starting or ending at each of the $n$ points
 $x_1 , \ldots , x_n$ is even.
\end {proposition}

The particular case where the set of points $\{ x_1, \ldots, x_n \}$ is the whole vertex set is again of some interest: The conditional distribution of the number of unoriented jumps on the edges given the occupation time field on the vertices is a collection of independent 
Poisson random variables with respective means $l(x_i) l(x_j)$, conditioned by the event that for all site $x$, the total number of jumps on the incoming edges at $x$ is even. 
This is exactly the random current distribution associated with the Ising model. For some further comments on this relation with random currents, the GFF and Ising, we refer to \cite {LuW}.

\subsection {Relation with Dynkin's isomorphism}

It should be of course noted that this decomposition is closely related Dynkin's isomorphism (see \cite {Dy,Dy2,MR} and the references therein), except that one here conditions here on the value of the square of the GFF instead of the value of the GFF itself. The previous result implies (when one only looks at occupation times and not at the loop-soup itself) 
that conditionally on the value of the square of the GFF at the set of points $\{ x_1, \ldots, x_n \}$, the square of the value of the GFF at the other points is the sum of the occupation times of the conditioned Poisson point process of excursions with an independent squared GFF in the remaining (smaller) domain.   

If one however conditions the GFF at the $n$ sites to be all equal to the same value $t$, then one can consider instead a graph where all these points are identified as a single point and note that
when the GFF on the new graph conditioned to have value $t$ at that point is distributed as the GFF on the initial graph, conditioned to have value $t$ at each of the $n$ points. 
One can apply the previous statement to that new graph and note that the conditioning on the event that the number of excursions-extremities  at each boundary site is even then disappears, because when there is just one such site, this number is anyway even (each excursion from this point to itself has two endpoints). Here it is however essential that the signs of all these values are the same (because if one identifies them into a single point, then they will anyway correspond to the same value of the GFF, not just to the same value of its square.

In summary, conditioning by the value of the square of the GFF gives rise to the parity conditioning, but it is also possible to condition on the actual 
value of the GFF and the parity conditioning becomes irrelevant when one looks at the occupation times only.
Note that Dynkin's isomorphism then follows, because in the latter case, the conditional distribution of the square of the GFF at the other points (which is therefore the square of the 
GFF in this smaller domain with boundary conditions given by these conditioned boundary values) will be the sum of the contribution of the 
loops that only visit those points (which is a squared GFF in the remaining domain) with the occupation time of the Poisson point process of excursions, while the conditioned GFF is a GFF with some prescribed boundary conditions, that can be viewed as the sum of a GFF in the complement of the set of marked points with the deterministic harmonic extension of these boundary values.

\subsection {The oriented case} 

One can follow almost word for word the same strategy to study the conditional distribution of oriented continuous-time loop-soups at $\alpha =1$ given their cumulated local time at sites. 
In that case, the excursions will be oriented, and the conditional distribution of the excursions away from these points will be a Poisson point process conditioned on the event that for each site, 
the number of incoming excursions is equal to the number of outgoing ones. 

The particular case where the set of points is the whole vertex set is again interesting. The conditional distribution of the set of jumps will be independent Poisson on each oriented edge, but conditioned 
on the fact that the number of incoming jumps at each site is going to its number of outgoing jumps. 
We leave all the details and further results to the interested reader.

\medbreak
\noindent
{\bf Note.} 
We found out that the recently posted preprint \cite {CL} describes some ideas that are similar to the present paper.
Our work was carried out totally independently of \cite {CL}.

\bigbreak
\noindent
{\bf Acknowledgements.}
The support and hospitality of SNF grant SNF-155922, NCCR-Swissmap, of the Clay foundation and the Isaac Newton Institute in Cambridge (where the present work has been carried out) are gratefully acknowledged.

\end{document}